\documentclass[11pt]{article}
\usepackage{latexsym}
\usepackage{amsfonts}
\makeatletter

\def\@footnotetext#1{\insert\footins{%
\footnotesize
    \interlinepenalty\interfootnotelinepenalty
    \splittopskip\footnotesep
    \splitmaxdepth \dp\strutbox \floatingpenalty \@MM
    \hsize\columnwidth \@parboxrestore
   \edef\@currentlabel{\csname p@footnote\endcsname\@thefnmark}\@makefntext
    {\rule{\z@}{\footnotesep}\ignorespaces
      #1\strut}}}
\def\abstract{\small\quotation{\hskip-\parindent\sc Abstract.}}

\def\classification{\@ifnextchar [{\@xfootnotenext}%
   {\begingroup\let\protect\noexpand
      \xdef\@thefnmark{}\endgroup
    \@footnotetext}}

\title {}

\begin{document}
\classification {{\it 2000 Mathematics
Subject Classification:} Primary 20E05, 20F05. }

\begin{center}
{\bf \Large   Automorphic orbits in free groups}

\bigskip

{\bf   Alexei G. Myasnikov \hskip 20pt Vladimir Shpilrain}

\end{center} 
\bigskip

\begin{abstract}
\noindent Let $F_n$ be the free group of a finite rank $n$. 
We study orbits $ Orb_{\phi}(u)$, where  $u$ is  an element of the  
group $F_n$, under the action of an automorphism $\phi$. If an orbit like that 
is finite, we
determine precisely  what its cardinality can  be  if $u$ runs through the 
whole group $F_n$, and $\phi$ runs through the whole  group  $Aut(F_n)$. 

 Another problem that we address here is related to Whitehead's 
algorithm that determines whether or not a given element of a free group 
of finite rank is an automorphic image of another given element. 
  It is
known that the first part of this algorithm (reducing a given free word
to a free word of minimum possible length by elementary Whitehead
automorphisms) is  fast (of quadratic time with respect to the
length of the word). On the other hand, the second part of the algorithm 
(applied to two words of the same minimum length) was always considered
very slow.
 We give here an improved algorithm for the  second part, and
we believe this algorithm always terminates in polynomial time
with respect to the length of the words. We prove that this is indeed
the case if the free group has rank 2.

\end{abstract}

\date{}

\bigskip

\noindent {\bf 1. Introduction }
\bigskip

  Let  $F_n$   be the free group of a finite rank  $n \ge 2$  with
a  set  $X = {\{}x_i {\}}, 1 \le i \le n$,  of free generators.
 Denote by  $ Orb_{\phi}(u)$ the orbit of an element $u$ of the free 
group $F_n$ under the action of an automorphism $\phi$. That 
is,  $ Orb_{\phi}(u)={\{}v \in F_n,  ~v= \phi^m(u)$ for some $m 
 \in {\bf Z_+}{\}}$.

One of the problems that  we address here is: how many elements can  a {\it finite}
orbit like that possibly have if $u$ runs through the whole group $F_n$,
and $\phi$ runs through the whole  group  $Aut(F_n)$? The answer is provided 
by the following theorem, in combination with a result of McCool \cite{McCool} 
(see also \cite{Khramtsov}): 
\medskip

\noindent {\bf Theorem 1.1.} In the free group $F_n$, there is
an orbit  $ Orb_{\phi}(u)$ of cardinality $k$ if and only if there is 
an element of order $k$ in the group $Aut(F_n)$. 
\medskip

 Thus, the question above is reduced to another question, of finding out 
what possible order can a torsion element of the group $Aut(F_n)$ have.
The latter was answered by McCool \cite{McCool}; more general results
were obtained later by Khramtsov \cite{Khramtsov}. We cite the relevant result 
in Section 2, after the proof of Theorem 1.1.  
\smallskip 

 It should be pointed out that the ``only if" part of our Theorem 1.1 is no longer valid if 
$\phi$ is  an arbitrary endomorphism. The following example is based on
the idea suggested by C. Sims. 
\smallskip 

\noindent {\bf Example.} In the free group $F_3$, let $\phi$ be the 
endomorphism that takes $x_1$ to $x_2^{-1} x_3$; ~$x_2$ to $x_1$;  
~$x_3$ to 1. Let $u=x_1x_2x_3$; ~then the cardinality of $ Orb_{\phi}(u)$ 
is 5, but there is no element of order 5 in the group $Aut(F_3)$.  

\medskip

 Another problem that we consider  here, is the following.
\smallskip 
 
 Let $u$ be an element of the free group $F_n$, whose length $|u|$ 
 cannot be decreased by any automorphism of $F_n$.  Let $A(u)$ denote the set 
 of elements ${\{}v \in F_n; |v| =  |u|, f(v)= u$  for some $f \in Aut(F_n){\}}$. 
How fast does  the cardinality of $A(u)$ grow as a function of $|u|$?
\smallskip 

 The set $A(u)$ is therefore an ``abridged" orbit $Orb_{Aut(F_n)}(u)$,
that includes only those automorphic images of $u$ that have the same length 
as $u$ does. 
\smallskip 

The  problem above was motivated by complexity issues for Whitehead's 
algorithm that determines whether or not a given element of a free group 
of finite rank is an automorphic image of another given element. It is 
known that the first part of this algorithm (reducing a given free word 
to a free word of minimal possible length by ``elementary" Whitehead 
automorphisms) is pretty fast (of quadratic time with respect to the 
length of the word). On the other hand, the second part of the algorithm 
 (applied to two words of the same minimum length) was always considered 
very slow. In fact, the procedure outlined in the original paper by 
Whitehead (see e.g. \cite{LS}), suggested this part of the algorithm to be of 
superexponential time with respect to the length of the words. 
 However, a standard trick in graph theory shows that there is an 
algorithm of at most exponential time (see Proposition 3.1 in Section 3). 
Moreover, in the case where the free group has rank 2, 
we were able to prove 
\medskip

\noindent {\bf Theorem 1.2.} Let $u \in F_2$ be a word whose length is 
irreducible by any automorphism of $F_2$ (in particular, $u$ is 
cyclically reduced). Then the  number of   automorphic images of $u$
that have the same length as $u$ does, is bounded by a polynomial function 
 of $|u|$. 
\medskip

 In fact, experimental data suggest that the number in the statement of Theorem 1.2 
has the (exact!) bound of  $8m^2-40m$ for $m \ge 9$, where $m= |u|$, 
but we were unable to prove that. 
\smallskip 

 Theorem 1.2 has the following 
\medskip

\noindent {\bf Corollary 1.3.} In the group  $F_2$,  Whitehead's
algorithm  terminates in  polynomial  time with
respect to the maximum length of the  two words in question. 
\medskip

 We do not know whether or not Theorem 1.2 and, therefore, 
Corollary 1.3 hold for free groups of bigger ranks. However, 
experimental data kindly provided by C. Sims allowed us to make
the following 

\medskip

\noindent {\bf Conjecture.} In the free group  $F_n$, the cardinality
of $A(u)$ is bounded by a polynomial of degree $2n-2$ in $|u|$, 
 provided the  length of $u$ is 
irreducible by any automorphism of $F_n$.  
\medskip 

  A most amazing thing is that, according to the experimental data
mentioned above,  the {\it maximum} cardinality of $A(u)$ that 
can actually occur under the irreducibility assumption in the  Conjecture, 
 appears to be {\it precisely}  a polynomial of degree $2n-2$ in $m=|u|$ for 
sufficiently large  $m$.  
 For  $n=2$, this polynomial, as we have already mentioned, is 
$8m^2-40m$ ~if ~$m \ge 9$.  For  $n=3$, the polynomial is  $48m^4 - 480m^3 +
1104m^2 - 672m$ if $m \ge 11$.  A particular element $u \in F_3$ of length 
$m$ whose orbit  $A(u)$ has the cardinality given by the latter polynomial, 
is, according to the same experimental data,  
$u = x_1^k x_2 x_1 x_2^{-1} x_1 x_2^2 x_3^2$, ~where $k = m-8$.  

\medskip 

 We also note that, in the case where the free group has rank 2 (but not in the 
general case),
the condition on $|u|$ to be irreducible by any automorphism 
can be relaxed to $u$ just being cyclically reduced. 
 If however we drop this latter condition, the
situation changes, and the number of automorphic images might 
become exponential: 
\medskip

\noindent {\bf Proposition 1.4.} The number of primitive elements
of length $m$
in the group $F_2$ (and therefore, in any group $F_n, n \ge 2$) is:
\smallskip

\noindent {\bf (a)} more than $\frac{8}{3\sqrt 3} \cdot (\sqrt 3)^m$ if $m$ 
is odd.

\noindent {\bf (b)} more than $\frac{4}{3} \cdot (\sqrt 3)^m$ if $m$ is
even.
\smallskip

\noindent {\bf (c)}  The number of {\it cyclically reduced}
primitive elements  of length $m \ge 1$ in the group $F_2$ is
$4 m \cdot \Phi(m)$, where $\Phi(m)$ is the Euler function of $m$,
i.e., the number of  positive  integers $<m$ relatively prime to $m$.
(Clearly, $\Phi(m)<m$). 
\medskip

Informally speaking, ``most" primitive elements  in $F_2$ 
are conjugates of primitive elements  of smaller length. 
This is not the case in $F_n$ for $n > 2$, where ``most"  primitive
elements are of the form $u \cdot  x_i^{\pm 1} \cdot v$ where  $u, v$ 
are   arbitrary elements that do not depend on $x_i$. 

 Proof of Proposition 1.4 is given in Section 4. \\

\noindent {\bf 2. Finite orbits}
\bigskip

We start with
\medskip

\noindent {\bf Proof of Theorem 1.1.}
\smallskip 

\noindent {\bf (1)} The ``only if" part  is a combination of an observation 
due to G.Levitt (see \cite{LevittLustig}) with a result of 
Bestvina and Handel \cite{BH}. Here is the  argument.
 Suppose that for 
some automorphism  $\varphi$ of the group  $F_n$, one has 
$\varphi^k(g)=g$  and $\varphi^q(g) \ne g$ for  $0 < q < k$.

Consider the action of  $\varphi$ on the subgroup $H = Fix
(\varphi^k)$    of all elements fixed by $\varphi^k$.  (This
subgroup is clearly invariant under  $\varphi$ since
$\varphi^k(\varphi(h))=\varphi(\varphi^k(h))=\varphi(h)$.)  Then 
$\varphi$  is  an  automorphism of $H$. Indeed, $\varphi$  is 
obviously surjective on  $H$ since for any $h \in H$, we have 
 $h =\varphi(\varphi^{k-1}(h))$.  
If $\varphi$  were
not injective on $H$, then we would have  $\varphi(h)=1$ for some 
$h \in H$, in which case $h$ could not be fixed by $\varphi^k$. 
   
Finally, $\varphi$ clearly  has order $k$ as an element
of the automorphism group $Aut(H)$. Since $H$ has rank at most $n$
by \cite{BH},  this  yields the ``only if" part of the theorem. $\Box$ 
\smallskip

\noindent {\bf (2)} To prove the ``if" part we need the following definition. 
A group $G$ satisfies
the {\it big powers } condition if for any tuple of elements $u_1,
\ldots, u_n$ from $G$ with $[u_i,u_{i+1}] \neq 1 ~(i = 1, \ldots, n-1),$ 
there is an integer $K$ such that for any integers $M_1, \ldots, M_n 
\geq K$, the following inequality holds:
$$u_1^{M_1}\ldots u_n^{M_n} \neq 1.$$
It is known that every free group
satisfies the big powers condition \cite{GB}. Now comes 
\smallskip

\noindent {\bf Lemma 2.1.} Let
$\phi$ be a nonidentical automorphism of $F_n.$ Then there exists an 
integer $K \ge 1$ such that for any $M_1, \ldots, M_n
\geq K$ the following inequality holds:
$$ \phi(x_1^{M_1} \ldots x_n^{M_n}) \neq x_1^{M_1} \ldots x_n^{M_n}.$$
\smallskip

{\bf Proof.} Suppose, by way of contradiction,  that for
any integer $K > 0$, there are integers $M_1(K), \ldots, M_n(K) \geq K$ such    
that    
$$ \phi(x_1^{M_1(K)} \ldots x_n^{M_n(K)}) = x_1^{M_1(K)} \ldots
x_n^{M_n(K)}. $$
It follows that
\begin{equation}
\label{eq:1}
\phi(x_1)^{M_1(K)}\ldots \phi(x_n)^{M_n(K)} x_n^{-M_n(K)} \ldots 
x_1^{-M_1(K)}=1       
\end{equation}
for all positive integers $K.$ As we have mentioned above, the free group 
$F_n$   
satisfies the big powers condition, therefore there are two commuting
consecutive factors in (1). Since $\phi$ is an automorphism, the 
only  consecutive factors which can possibly commute are $\phi(x_n)^{M_n(K)}$ and
$x_n^{M_n(K)}.$ It follows that $\phi(x_n) = x_n$ and (1) 
takes the form    

$$\phi(x_1)^{M_1(K)}\ldots \phi(x_{n-1})^{M_{n-1}(K)}
x_{n-1}^{-M_{n-1}(K)} \ldots x_1^{-M_1(K)}= 1.$$

Upon repeating the argument above, we get $\phi(x_i) = x_i$ for all $i = 1,  
\ldots,n,$ i.e., $\phi$ is identical. 
This contradiction proves the lemma. $\Box$
\smallskip

 We now continue with our proof of the ``if" part. Given $k > 1$
and an automorphism $\varphi$ of order $k$  of the group $F_n$, we
are going to  find an element $u \in F_n$,  so that the 
orbit $ Orb_{\varphi}(u)$ has cardinality $k$.

 If $\varphi$ is a permutation on the set $\{x_1^{\pm 1},...,
x_n^{\pm 1}\}$, then any element of the form  $u=x_1^{M_1} \cdot ... \cdot
x_n^{M_n}, ~M_i \ne 0,$ would do.  If not, then there is at least one free
generator, say,  $x_1$, such
  that  $\varphi(x_1)$ has length at least 2.  Let $u=
x_1^{M_1} \cdot ...  \cdot x_n^{M_n}$.
Then, by Lemma 2.1, for some choice of $K \ge 1$, for any 
 $M_1, \ldots, M_n \geq K$ we have $\varphi(u) \ne u$.

 Similarly, for any $m, ~1<m<k$, we can construct an element $u_m$   such 
  that $\varphi^m(u_m) \ne u_m$. Every $u_m, ~m \ge 2,$ is chosen to be of the 
form $u_m=x_1^{M_{1,m}} \cdot ...  \cdot x_n^{M_{n,m}}$ with $min_i M_{i,m} >
max_i M_{i,m-1}$, ~and $\varphi^m(u_m) \ne u_m$ (the latter is possible by Lemma 2.1).

 Obviously, with this choice of $M_{i,j}$ we will also have 
$\varphi^j(u_j) \ne u_j$ for any $j \le m$. Therefore, for $u=u_k$, 
the orbit $ Orb_{\varphi}(u)$ will  have cardinality $k$. $\Box$ 
        
\medskip

 We note that possible values of the  order  of a torsion
element of the group $Aut(F_n)$ are described, according to 
\cite{McCool} and \cite{Khramtsov}, as follows. 
Pick a positive integer $k=p_1^{\alpha_1} \cdot ... \cdot 
p_s^{\alpha_s}$, where $p_1, ..., p_s$ are different primes. 
There is an element of order $k$ in the group $Aut(F_n)$ if 
and only if $\sum_{i=1}^{s} (p_i^{\alpha_i}-p_i^{\alpha_{i-1}}) 
\le n$. For example, if $k=15=3 \cdot 5$; then the sum above 
becomes (3-1)+(5-1) =6. Therefore, there is an automorphism 
of order 15 in the group $Aut(F_n)$ for $n \ge 6$, but not for
$n \le 5$.

 We also note that Levitt and  Nicolas \cite{Levitt}  proved 
that the {\it maximum} order (call it $H(n)$) of a torsion element of  $Aut(F_n)$
is the same as that of a torsion element of  $GL_n({\bf Z})$, 
with the exception of $n$ =2, 6, and 12. They also established 
the asymptotic of this  function  by showing 
$log H(n) \sim \sqrt{n \cdot log ~n}$.  \\

\noindent {\bf 3. Whitehead's algorithm revised}
\bigskip

 In this section, we study complexity of Whitehead's
algorithm that determines whether or not a given element of a free
group  of finite rank is an automorphic image of another given
element.

 It is
known that the first part of this algorithm (reducing a given free word
to a free word of minimum possible length by ``elementary" Whitehead
automorphisms) is pretty fast (of quadratic time with respect to the
length of the word). On the other hand, the second part of the algorithm
(applied to two words of the same minimum length) was always considered
very slow. In fact, the procedure outlined in the original paper by
Whitehead \cite{W}, suggested this part of the algorithm to be of
superexponential time with respect to the length of the words.
Indeed, given a word $u$, the procedure calls for constructing a
graph whose  vertices correspond to all words of length $|u|$. 
That means, the number of vertices is an exponential function of 
$|u|$. After that, for every vertex of the graph,  one 
constructs edges incident to this vertex as follows: an edge
connects   this particular vertex to  another vertex if and only 
if there is an elementary Whitehead automorphism that takes
one of the corresponding words to the other. Finally, to find out
if there is an automorphism that takes the word $u$ to 
another given word $v$ of the same length, one has to check 
all the paths in the graph that start at the vertex that 
corresponds to $u$, and see if some of them leads to the vertex that 
corresponds to $u$. The number of paths in  a graph is, in 
general, an exponential function of the number of vertices, 
therefore this algorithm is, in general, of superexponential complexity with 
respect to the length of the word $u$.
\smallskip

 It is possible however to skip some steps in this algorithm 
and get the following 
\medskip

\noindent {\bf Proposition 3.1.} Let $N$ be the number of 
automorphic images of $u \in F_n$ that have the same length as $u$
does.  Then, given an element $v$ of length $|u|$, one can 
 decide in  linear  time with respect to $N$, whether or not 
$v$ is an automorphic image of $u$. 
\medskip

\noindent {\bf Proof.} We are going to use the {\it backtracking} 
method which is a well-known procedure in graph theory for
searching a  tree. 

 Starting with the vertex that 
corresponds to $u=u_0$, we are building a tree as follows. 
(We use the same notation for words and corresponding vertices 
when there is no ambiguity). 
\smallskip

\noindent {\bf (1)} Apply  an arbitrary elementary Whitehead 
automorphism to $u_0$; if a new word $u_1$ of the same length 
is obtained, plot the corresponding vertex and connect it to 
$u_0$. If not, then apply another elementary Whitehead 
automorphism, until you get a new word $u_1$ of the same length. 
(Note that the total number of those automorphisms $C=C(n)$ is
finite and depends on the rank $n$ of the group $F_n$ only). 
\smallskip

\noindent {\bf (2)} Continue the same process. That is, suppose we
have obtained  a  word  $u_i, ~i>0$, at the previous
step. This time ``a new word"
would mean a word different  from all the words obtained at 
previous steps. 

 If none of the  elementary Whitehead automorphisms produces 
a new word, then do ``backtracking", i.e., return to the 
 word obtained at the immediately preceding step, and 
repeat the same process. 
\smallskip

 In the end (i.e., when no new word can be obtained from any 
of the ``old" words), we shall obviously have a {\it spanning tree} 
of the graph described before the statement of Proposition 3.1.
It will therefore have $N$ vertices and $N-1$ edges. Furthermore, 
in the course of constructing this tree, we did not   traverse  
  any of the edges more than twice (once in each direction). 

 Thus, the time we need to construct this tree, is no more than 
$C \cdot N$, where $C$ is the constant mentioned above. 
 Once the tree is constructed, it will take just $N$ more steps 
to find out if the vertex corresponding to the word $v$ is 
among the vertices. Or, we can perform the check every time 
we get a new vertex, because once we get $v$, we can stop. $\Box$
 
\medskip

 Thus, the speed of Whitehead's algorithm is determined by 
the number  of automorphic images of an element $u \in F_n$ that
have the same length as $u$ does. Therefore, Theorem 1.2 will 
imply that, in the case where the free group has rank 2, 
Whitehead's algorithm does, in fact, terminate in  polynomial 
time with respect to the length of the words in question. 

 We are now ready for 

\medskip

\noindent {\bf Proof of Theorem 1.2.} Throughout the proof, 
we shall call ``length-preserving" those automorphisms of $F_2$
that are permutations on the set $\{x, x^{-1}, y, y^{-1}\}$. 
There are 8 of them, so whenever we count the number of 
automorphic images of a particular element ``up to a
length-preserving automorphism", it means the upper bound for
such a number should be multiplied by 8.

 Let $M=|u|$. Let $k$ be the sum of exponents on  $x$ 
in the word $u$, and $l$  the sum of exponents on  $y$.
Upon applying a length-preserving automorphism  
if necessary, we may assume that $k, l \ge 0$. First, we are going to 
establish the result of Theorem 1.2 for $u \notin [F_2,F_2]$, so we assume
that $k, l$ are not both 0. In this case, the result 
will follow from the following observations.
\smallskip

\noindent {\bf (1)}
 For a word of length $M$, there are $\sum_{i=0}^M (i+1) =
 \frac{1}{2} (M+1)(M+2)$ possible pairs $(k, l)$ with $k, l \ge
0; ~k + l \le M$. 
\smallskip

\noindent {\bf (2)} It is well known  (see e.g. \cite{MKS}) that  
the group $Aut(F_2)$ is generated by inner automorphisms, by 
 3 length-preserving  automorphisms $\pi : x \to y, ~y \to x$; 
 $~\sigma_x : x \to x^{-1}, ~y \to y$;  $~\sigma_y : x \to x, ~y \to y^{-1}$, 
  and by the following two: $\alpha : x \to xy, ~y \to y$, and 
~$\beta :  x \to x, ~y \to yx$.  The subgroup $H$ of $Aut(F_2)$  generated by  
    $\alpha$ and ~$\beta$ can be mapped onto  
  $SL_2({\bf Z})$. Under this epimorphism, $\alpha$ and ~$\beta$ 
correspond to the matrices 
$\left(\begin{array}{cc} 1 & 0  \\ 1 & 1 
 \end{array}\right)$ and $\left(\begin{array}{cc} 1 & 1  \\ 0 & 1 
 \end{array}\right)$, respectively. The kernel of this epimorphism 
 is generated  (as a normal subgroup)  by the inner automorphism 
induced by the element $[x,y]$; in particular, every automorphism 
in the kernel is inner.

 Furthermore, relations between generators of $Aut(F_k)$ given in 
\cite[Section 3.5, Theorem N1]{MKS} show that in any product of 
automorphisms $\alpha^{\pm 1}$, ~$\beta^{\pm 1}$, ~$\pi$, $~\sigma_x$, and 
$~\sigma_y$,  ~automorphisms $\alpha^{\pm 1}$ and  ~$\beta^{\pm 1}$ 
can be collected on the right.  This, together with the fact that 
 the subgroup of inner automorphisms of $F_2$ is normal in $Aut(F_2)$, 
 implies that applying an automorphism of $F_2$ amounts to first 
applying an automorphism from the subgroup $H$   generated by  
    $\alpha$ and ~$\beta$, then  a length-preserving  automorphism, 
 and, finally, an inner automorphism. 

Therefore, to bound the  number of
cyclically reduced automorphic images of $u$ with the same non-zero  
vector  $(k, l)$ of exponent sums, it is sufficient to bound the  
number of matrices from $SL_2({\bf Z})$ that fix the vector $(k, l)$ 
acted upon by right multiplication,  and then multiply this number 
by $M$ (the number of cyclic permutations of a word of length 
 $M$). Furthermore, up to a length-preserving automorphism, 
every automorphism from the group  $H$  corresponds to a matrix from $SL_2({\bf
Z})$ whose elements in the first row are of different signs, say,
the element in the upper left corner is non-negative,  and the element in the 
upper right corner is non-positive. 
(Elements in the first row correspond to the image of $x$).

\smallskip

\noindent {\bf (3)} Thus, what is left to do now is to count the number 
of matrices in $SL_2({\bf Z})$ whose elements in the first row  
are  of different signs, that fix a given non-zero vector $(k, l)$ with 
$k, l \ge 0$. The  computation here is straightforward.  Let $A= 
\left(\begin{array}{cc} a_{11} & a_{12}  \\ a_{21} & a_{22} 
 \end{array}\right)$ 
be a matrix from $SL_2({\bf Z})$ with $a_{11}  \ge 0, a_{12} \le 0$, 
which fixes  a vector $(k, l)$. Then we have the following 
system  of equations in $a_{ij}$:

\noindent $k \cdot a_{11} + l \cdot a_{21} = k$; 
~$k \cdot a_{12} + l \cdot a_{22} = l$; 
~$a_{11} a_{22} - a_{12} a_{21} = 1$.

Suppose first that both $k, l \ne 0$.  
Then from the first equation we get $a_{21} =\frac{k}{l} - \frac{k}{l}\cdot  
a_{11}$, ~and from the second equation $a_{22} = 1- \frac{k}{l} \cdot 
a_{12}$. Plug this into the third equation and simplify: 
$l \cdot a_{11} - k \cdot a_{12}= l$. Since $k, l >0, ~ a_{11}  \ge 0, a_{12} \le 0$, 
this gives either $a_{12}= 0, a_{11} =1$,  
or $a_{11}= 0, ~a_{12}= -\frac{l}{k} $. In the former case, we get 
$a_{22}=1, a_{21}= 0$.
In the latter case, $a_{21}= \frac{k}{l}, a_{11} =0, a_{22}= 2$. 

 Now suppose, say,  $k =0$. Then $a_{21}= 0, a_{22}= 1,  a_{11}=1$, 
whereas  $a_{12}$ can be arbitrary. However, we can show that, 
should the automorphism corresponding to the matrix $A$ preserve 
the length of $u$, the absolute value of $a_{12}$ cannot be greater 
than $2|u|$. Indeed, let $K=a_{12}$; then the automorphism corresponding 
to the matrix $A$ is $\alpha^K$, i.e., it takes $x$ to $xy^K$, ~$y$ to $y$. 
Suppose $K>2|u|$; we may assume that $u$ has at least one occurence of 
$x$. Then  $\alpha^K(u)$ has a subword $xy^K$ (before cancellation). 
Since we have assumed that $\alpha^K(u)$ has the same length as $u$ does, 
more than half of  $y^K$ should cancel out. This implies that, 
in the word $u$ itself, there is a subword $y^{-N}$ with $N \ge \frac{K+1}{2}$. 
This is a contradiction since $K>2|u|$.  

\smallskip

 Thus, in any of the considered cases, we have no more than  $2|u|$   
different matrices from $SL_2({\bf Z})$ that fix a given non-zero 
vector $(k, l)$. 
\smallskip

 Summarizing the observations (1), (2), (3), we see that 
the number of cyclically reduced automorphic images of $u$ 
of length $M=|u|$ is no more than $c \cdot M^4$ for some 
constant $c$ independent of $u$. This completes the proof in the case where 
$u \notin [F_2,F_2]$. 
\smallskip

 Now let $u \in [F_2,F_2]$. In this case, we are going to use induction
on the length of $u$. To make the induction work, we are going 
to prove the following somewhat stronger claim. 
\medskip

\noindent {\bf Proposition 3.2.} Let $u \in [F_2, F_2]$ be cyclically 
reduced. For any positive 
integer $K$,  the number of elements $v \in F_2$  such 
  that $v = \phi(u)$ for some     $\phi \in Aut(F_2)$ and ~$|v| = |u| +K$, 
is less than $c \cdot 3^K \cdot  (|u|+K)^4$ for some constant $c$ independent of 
$u$ and ~$K$.
\medskip 

\noindent {\bf Proof.} The basis of induction $u=[x,y]$ is almost obvious. This element
is fixed by any automorphism from $H$ (recall that $H$ is the subgroup of $Aut(F_2)$ 
generated by two automorphisms, 
$\alpha : x \to xy, ~y \to y$, and ~$\beta :  x \to x, ~y \to yx$),  and therefore,
to count the number of elements $v \in F_2$  such 
  that $v = \phi(u)$ for some     $\phi \in Aut(F_2)$ and ~$|v| = |u| +K$, we just
have to count (up to a length-preserving automorphism) the number of conjugates of $u$ 
of length up to $|u| +K$. This latter 
number is no bigger than the number of different elements of length $[K/2]$ in the group $F_2$, 
i.e., equals $3^{[K/2]}$.

\smallskip

 For the induction step, we first assume that $u$ has a subword of the form 
$[x^{\pm 1}, y^{\pm 1}]$. Then, upon applying a length-preserving automorphism
if necessary, we may assume that $u$ has a subword $[x,y]$. Then a cyclic permutation
of $u$ has the form $[x,y] w$, with no cancellation between 
$[x,y]$ and $w$. Thus, by the remarks in the beginning of the proof of Theorem 1.2, 
we may assume that $u=[x,y] w$. Let $\phi$ be an arbitrary automorphism from $H$.
Recall that every automorphism in $H$ fixes $[x,y]$. 

 Assume first that $w$ is cyclically reduced. 
 We have two possibilities:
\smallskip

  \noindent {\bf (1)} $[x,y]$ is entirely canceled out by   $\phi(w)$.  
Then, since $\phi(u)=  [x,y] \phi(w)$, we see that, if
 $|\phi(u)| = |u| +K$,  we must have  $|\phi(w)| = |w|+K+4$. 
 By the inductive assumption, the number of automorphic images of $w$ with this property 
is no more than $c \cdot 3^{K+4} \cdot |w+K+4|^4$ for some constant $c$ independent 
of $w$ and ~$K$. Similar result for $u$ now follows.  

\smallskip

  \noindent {\bf (2)} Only part of $[x,y]$ cancels out (this includes 
the case where nothing cancels out).  Then, since $\phi(u)= 
  [x,y] \phi(w)$ ~and since an element of the commutator subgroup must have 
an even length, we see that, if
 $|\phi(u)| = |u| +K$,  then either $|\phi(w)| = |w|+K+2$, or $|\phi(w)| = |w|+K$. 
 By the inductive assumption, the number of automorphic images of $w$ with this property 
is no more than $c \cdot 3^{K+2} \cdot |w+K+2|^4$ (respectively,  
$c \cdot 3^K \cdot |w+K|^4$) for some constant $c$ independent of $w$ and ~$K$.
 Similar result for $u$ now follows.  

\smallskip

 If $w$ is not cyclically reduced, i.e., if $u=[x,y] g w' g^{- 1}$, then we consider
a cyclic permutation of $u$: ~$u'= g^{- 1}[x,y] g w'= [x^{g^{- 1}}, y^{g^{- 1}}] w'$, 
where we can assume $w'$ to be cyclically reduced. Now we apply essentially the same
argument to $u'$ as we have just applied to $u$, upon replacing the subgroup $H$ 
of automorphisms by the left coset $i_{g} H$, where $i_{g}$ is the inner 
automorphism induced by the element $g$. (Applying an automorphism from $i_{g} H$ 
is equivalent to first applying conjugation by $g$, and then applying     
an automorphism from $H$.) 

 Since the group of  inner automorphisms is normal in $Aut(F_2)$, observation (2)
in the beginning of the proof of Theorem 1.2 remains  
 valid upon replacing $H$ by $i_{g} H$. That is, every  automorphism from $Aut(F_2)$
 is a product of an automorphism from the  coset $i_{g} H$ and an inner automorphism.    
Since every automorphism from $i_{g} H$ fixes the element 
$[x^{g^{- 1}}, y^{g^{- 1}}]$, the same argument as above completes the 
proof in this case.  

\smallskip

 Suppose now that $u$ does not have a subword of the form 
$[x^{\pm 1}, y^{\pm 1}]$, but does  have a subword of the form 
$x^{\pm 1} y^{\pm 1} x^{\mp 1}$. Then, upon applying a length-preserving automorphism
if necessary, we may assume that $u$ has a subword $x y x^{- 1}$. Thus, a cyclic permutation
of $u$ has the form $x y x^{- 1} w$, with no cancellation. Then we can write $u$ as 
$u=[x,y] y w$. Note that the word $y w$ has smaller length than $u$ does, 
 and we can assume that $y w$ is cyclically reduced, for if it was not, $w$ would 
end with $y^{-1}$, and then a cyclic permutation
of $u$ would be of the form $y^{-1}x y x^{- 1} w'=[y^{-1}, x] w'$, and therefore this case
 would be reduced to the previous one. 

 Thus, we can apply the inductive assumption to this word $w$, and  the same 
argument as above will work in this case as well.  

\smallskip

 Finally, suppose that $u$ does not have a subword of the form 
$x^{\pm 1} y^{\pm 1} x^{\mp 1}$. Then $u$ must have a subword of the form 
$x^{\pm 1} y^k x^{\mp 1}$ ~for some $k \ne 0, \pm 1$. We can assume, upon 
applying a length-preserving automorphism and a cyclic permutation if necessary,
that $u=x y^k x^{- 1} w, ~k > 1$.   Then we can write $u=[x, y] y x y^{k- 1} x^{- 1} w$. 
 Now the word $y x y^{k- 1} x^{- 1} w$ has the same length as $u$ does, but it 
has the subword $x y^{k-1} x^{- 1}$. Also, we can assume that $y w$ is cyclically 
reduced, for if it was not, $w$ would 
end with $y^{-1}$, and then a cyclic permutation
of $u$ would be of the form $y^{-1}x y^k x^{- 1} w'$, i.e., it would begin with 
$y^{-1}x y$,  and therefore this case
 would be reduced to   one of the previously considered. 
 An obvious inductive argument  now completes the proof. $\Box$ \\

\noindent {\bf 4. Primitive elements of $F_2$}
\bigskip

 In this section, we give bounds for the total number of primitive elements
of a given length $m$ in the group $F_2$, and a precise number of  {\it
cyclically reduced}  primitive elements of length $m$. (Note that the total 
number of
elements of length $m$ in the group $F_2$ is $\frac{4}{3} \cdot 3^m$).
\medskip

\noindent {\bf Proof of Proposition 1.4.} Let $x$ and  $y$ be generators of
$F_2$.

\noindent {\bf (a)} Suppose $m$ is odd. Then  any conjugate of $x^{\pm 1}$,
as well as of $y^{\pm 1}$, by an element of length $k= (m-1)/2$, is a
primitive element of length $m$ (assuming there are no cancellations in the middle).
The number of  elements like that  in the group $F_2$ is $2 \cdot 3^{k-1}$, whence
the result.
\medskip

\noindent {\bf (b)} If $m$ is even, then counting conjugates of $x^{\pm
1}y$  and  $x y^{\pm 1}$ by elements of length $(m-2)/2$ yields the result.
\medskip

\noindent {\bf (c)} The result of this part will follow from a well-known 
fact about primitive elements of $F_2$ (see \cite{CMZ} or \cite{OZ}):

\noindent -- for any pair $\{k,l\}$ of integers with $(k,l)= 1$, there
is exactly one cyclically reduced  primitive element of $F_2$ whose exponent sum
on $x$ is $k$ and the  exponent sum on $y$ is $l$.

 Thus, the number of cyclically reduced  primitive elements of $F_2$ of 
length
$m$ is 8 times the number of pairs  $\{k,l\}$ of {\it positive} integers
with $(k,l)= 1$, $~k <l, ~k+l= m$. The latter number is obviously equal to
$\frac{1}{2}\Phi(m)$, where $\Phi(m)$ is the number of positive integers $<m$ relatively prime to $m$.
$\Box$\\

\noindent {\bf Acknowledgement} 
\medskip 

 We are  indebted to A. D. Myasnikov and  C. Sims for providing us with 
 useful experimental data.

\medskip
\noindent 
 Department of Mathematics, The City  College  of New York, New York, 
NY 10031 

\medskip

\noindent {\it e-mail addresses\/}:  alexeim@att.net,
 shpil@groups.sci.ccny.cuny.edu 
\medskip

\noindent {\it http://www.grouptheory.info   }

\end{document}